\documentclass{amsart}
\usepackage{tipa}
\usepackage{amssymb}
\usepackage{stmaryrd}
\usepackage{graphicx}
\usepackage{amsmath,amsfonts}
\usepackage{mathtools}

\newtheorem{theorem}{Theorem}[section]

\theoremstyle{definition}

\newtheorem{example}[theorem]{Example}

\numberwithin{equation}{section}

\begin{document}

\title[Heron triangle and rhombus pairs]{Heron triangle and rhombus pairs with a common area and a common perimeter}

\author{Yong Zhang}
\address{School of Mathematics and Statistics, Changsha University of Science and Technology,
Changsha 410114, People's Republic of China}
 \email{zhangyongzju$@$163.com}

\author{Junyao Peng}
\address{School of Mathematics and Statistics, Changsha University of Science and Technology,
Changsha 410114, People's Republic of China}
 \email{junyaopeng906$@$163.com}

\thanks{This research was supported by the National Natural Science Foundation of China (Grant No.~11501052).}

\date{}


\begin{abstract}
By Fermat's method, we show that there are infinitely many Heron
triangle and $\theta$-integral rhombus pairs with a common area and
a common perimeter. Moreover, we prove that there does not exist any
integral isosceles triangle and $\theta$-integral rhombus pairs with
a common area and a common perimeter.
\end{abstract}

\maketitle

\section {Introduction}
We say that a Heron (resp. rational) triangle is a triangle with
integral (resp. rational) sides and integral (resp. rational) area.
And a rhombus is $\theta$-integral (resp. $\theta$-rational) if it
has integral (resp. rational) sides, and both $\sin\theta$ and
$\cos\theta$ are rational numbers.

In 1995, R. K. Guy \cite{Guy} introduced a problem of Bill Sands,
that asked for examples of an integral right triangle and an
integral rectangle with a common area and a common perimeter, but
there are no non-degenerate such. In the same paper, R. K. Guy
showed that there are infinitely many such integral isosceles
triangle and rectangle pairs. In 2006, A. Bremner and R. K. Guy
\cite{Bremner-Guy} proved that there are infinitely many such Heron
triangle and rectangle pairs. In 2016, Y. Zhang \cite{Zhang} proved
that there are infinitely many integral right triangle and
parallelogram pairs with a common area and a common perimeter. At
the same year, S. Chern \cite{Chern} proved that there are
infinitely many integral right triangle and $\theta$-integral
rhombus pairs. In a recent paper, P. Das, A. Juyal and D. Moody
\cite{Das-Juyal-Moody} proved that there are infinitely many
integral isosceles triangle-parallelogram and Heron triangle-rhombus
pairs with a common area and a common perimeter.

By Fermat's method \cite[p. 639]{Dickson2}, we can give a simple
proof of the following result, which is a corollary of Theorem 2.1
in \cite{Das-Juyal-Moody}.
\begin{theorem}
There are infinitely many Heron triangle and $\theta$-integral
rhombus pairs with a common area and a common perimeter.
\end{theorem}

But for integral isosceles triangle and $\theta$-integral rhombus
pair, we have
\begin{theorem}
There does not exist any integral isosceles triangle and
$\theta$-integral rhombus pairs with a common area and a common
perimeter.
\end{theorem}

\section {Proofs of the theorems}
\begin{proof}[\textbf{Proof of Theorem 1.1.}] Suppose that the Heron triangle has sides $(a,b,c)$, and the $\theta$-integral
rhombus has side $p$ and intersection angle $\theta$ with
$0<\theta\leq \pi/2$. By Brahmagupta's formula, all Heron triangles
have sides
\[(a,b,c)=((v+w)(u^2-vw),v(u^2+w^2),w(u^2+v^2)),\]
for positive integers $u, v, w,$ where $u^2>vw$.

Noting that the homogeneity of these sides, we can set $w=1$, and
$u,v,p$ be positive rational numbers. Now we only need to study the
rational triangle and $\theta$-rational rhombus pairs with a common
area and a common perimeter, then we have
\begin{equation}\label{Eq21}
\begin{cases}
\begin{split}
&uv(v+1)(u^2-v)=p^2\sin{\theta},\\
&2u^2(v+1)=4p.
\end{split}
\end{cases}
\end{equation}
Since both $\sin\theta$ and $\cos\theta$ are rational numbers, we
may set
\[\sin{\theta}=\frac{2t}{t^2+1},\cos{\theta}=\frac{t^2-1}{t^2+1},\]
where $t \geq1$ is a rational number. For $t=1$, $\theta=\pi/2$,
this is the case studied by R. K. Guy \cite{Guy}. Thus it only needs
to consider the case $t>1$.

Eliminating $p$ in Eq. (\ref{Eq21}), we have
\[\frac{u(v+1)(2t^2u^2v-tu^3v-2t^2v^2-tu^3+2u^2v-2v^2)}{2(t^2+1)}=0.\]
Let us study the rational solutions of the following equation
\begin{equation}\label{Eq22}
2t^2u^2v-tu^3v-2t^2v^2-tu^3+2u^2v-2v^2=0,
\end{equation}
and solve it for $v$, we get
\[v=\frac{(2t^2u-tu^2+2u\pm\sqrt{g(t)})u}{4(t^2+1)},\]
where
\[g(t)=4u^2t^4-4u(u^2+2)t^3+u^2(u^2+8)t^2-4u(u^2+2)t+4u^2.\]
In view of $v$ is a positive rational number, then $g(t)$ should be
a rational perfect square. So we need to consider the rational
points on the curve
 \[\mathcal{C}_1:~s^2=g(t).\]
The curve $\mathcal{C}_1$ is a quartic curve with a rational point
$P=(0,2u)$. By Fermat's method \cite[p. 639]{Dickson2}, using the
point $P$ we can produce another point $P'=(t_1,s_1)$, which
satisfies the condition $t_1s_1\neq0$. In order to construct a such
point $P'$, we put
\[s= rt^2+ qt +2u,\]
where $r,q$ are indeterminate variables. Then
\[s^2-g(t)=\sum_{i=1}^4A_it^i,\]
where the quantities $A_i = A_i(r,q)$ are given by
\[\begin{split}
A_1=&4u^3+4qu+8u,\\
A_2=&-4u^2+4ru+q^2-8u^2,\\
A_3=&4u^3+2rq+8u,\\
A_4=&r^2-4u^2.
\end{split}\]
The system of equations $A_3=A_4=0$ in $r, q$ has a solution given
by\[r=-2u,q=u^2+2.\] This implies that the equation
\[s^2-g(t)=\sum_{i=1}^4A_it^i=0\] has the rational roots $t=0$ and
\[t=\frac{2u(u^2+2)}{3u^2-1}.\]
Then we have the point $P'=(t_1,s_1)$ with
\[\begin{split}
t_1=&\frac{2u(u^2+2)}{3u^2-1},\\
s_1=&-\frac{2u(u^6-4u^4+14u^2+3)}{(3u^2-1)^2}.
\end{split}\]
Put $t_1$ into Eq. (\ref{Eq22}), we get
\[v=\frac{u^2(u^2+2)}{4u^4+1}.\]
Hence, the rational triangle has rational sides $(a,b,c)=$
\[\bigg(\frac{u^2(3u^2-1)(u^4+6u^2+1)}{(4u^2+1)^2}, \frac{u^2(u^2+2)(u^2+1)}{4u^2+1},\frac{u^2(u^6+20u^4+12u^2+1)}{(4u^2+1)^2}\bigg).\]
From the equation $2u^2(v+1)=4p$ and $\sin\theta=\frac{2t}{t^2+1}$,
we obtain the corresponding rhombus with side
\[p=\frac{(u^4+6u^2+1)u^2}{2(4u^2+1)},\] and the intersection angle
\[\theta=\arcsin\frac{4u(u^2+2)(3u^2-1)}{(4u^2+1)(u^4+6u^2+1)}.\]

In view of $u,v, p$ are positive rational numbers, $0<\sin\theta<1$,
$u^2>v$ and $t_1>1$, we get the condition
\[u>\frac{\sqrt{3}}{3}.\] Then for positive rational number
$u>\frac{\sqrt{3}}{3}$, there are infinitely many rational triangle
and $\theta$-rational rhombus pairs with a common area and a common
perimeter. Therefore, there are infinitely many such Heron triangle
and $\theta$-integral rhombus pairs.
\end{proof}

\begin{example} (1) If $u=1,$ we have a Heron triangle with
sides $(8,15,17)$, and a $\theta$-integral rhombus with side $10$
and the smaller intersection angle $\arcsin(3/5),$ which have a
common area $60$ and a common perimeter $40.$

(2) If $u=2,$ we have a Heron triangle with sides $(1804,2040,1732)$
and a $\theta$-integral rhombus with side $1394$ and the smaller
intersection angle $\arcsin(528/697),$ which have a common area
$1472064$ and a common perimeter $5576.$
\end{example}

\begin{proof}[\textbf{Proof of Theorem 1.2.}] As before, we only
need to consider the rational isosceles triangle and
$\theta$-rational rhombus pairs. As in \cite{Das-Juyal-Moody}, we
may take the equal legs of the isosceles triangle to have length
$u^2+v^2$, with the base being $2(u^2-v^2)$ and the altitude $2uv$,
for some rational $u, v$. The area of the isosceles triangle is
$2uv(u^2-v^2)$, with an perimeter of $4u^2$.

Let $p$ be the length of the side of the rhombus, and $\theta$ its
smallest interior angle. For $\theta$-rational rhombus, we have the
perimeter $4p$ and area $p^2\sin\theta$, where
$\sin\theta=2t/(t^2+1)$, for some $t\geq1$.

If the rational isosceles triangle and $\theta$-rational rhombus
have the same area and perimeter, then
\begin{equation}\label{Eq23}
\begin{cases}
\begin{split}
&2uv(u^2-v^2)=p^2\sin\theta,\\
&4u^2=4p.
\end{split}
\end{cases}
\end{equation}
From Eq. (\ref{Eq23}), we obtain
\[\frac{2u(v(u-v)(u+v)t^2-u^3t+v(u-v)(u+v))}{t^2+1}=0.\]
It only needs to consider $v(u-v)(u+v)t^2-u^3t+v(u-v)(u+v)=0.$ If
this quadratic equation has rational solutions $t$, then its
discriminant should be a rational perfect square, i.e.,
\[u^6-4u^4v^2+8u^2v^4-4v^6=w^2.\]
Let $U=u/v,W=w/v^3$, we have
\[W^2=U^6-4U^4+8U^2-4.\]
This is a hyperelliptic sextic curve of genus 2. The rank of the
Jacobian variety is 1, and Magma's Chabauty routines determine the
only finite rational points are
\[(U,W)=(\pm1;\pm1),\] which
lead to \[(u,w)=(\pm v,\pm v^3),\] then we get
\[u^3t=0.\]
So Eq. (\ref{Eq23}) does not have nonzero rational solutions, which
means that there does not exist any integral isosceles triangle and
$\theta$-integral rhombus pairs with a common area and a common
perimeter.
\end{proof}

\vskip20pt
\bibliographystyle{amsplain}

\end{document}